\documentclass[10pt]{article}

\usepackage{amscd,amsmath, amssymb, fancyhdr, url, graphicx}

\usepackage[h=2.3em]{diagrams}

\numberwithin{equation}{section}

\newcommand{\version}{version 0.9,\ \ Apr 13, 2016}

\makeatletter
\def\x@arrow{\DOTSB\Relbar}
\def\xlongrightarrowfill@{\arrowfill@\relbar\relbar\longrightarrow}
\newcommand{\xlongrightarrow}[2][]{%
        \ext@arrow 0099\xlongrightarrowfill@{#1}{#2}}
\makeatother

\def\eqref#1{(\ref{#1})}
\newcommand{\goth}{\mathfrak}
\newcommand{\g}{{\mathfrak g}}
\newcommand{\arrow}{{\:\longrightarrow\:}}
\newcommand{\Z}{{\Bbb Z}}
\def\C{{\Bbb C}}

\newcommand{\R}{{\Bbb R}}
\newcommand{\Q}{{\Bbb Q}}
\renewcommand{\H}{{\Bbb H}}

\def\1{\sqrt{-1}\:}

\newcommand{\cntrct}                
{\hspace{2pt}\raisebox{1pt}{\text{$\lrcorner$}}\hspace{2pt}}

\renewcommand{\phi}{\varphi}
\renewcommand{\epsilon}{\varepsilon}
\renewcommand{\geq}{\geqslant}
\renewcommand{\leq}{\leqslant}

\newcommand{\Teich}{\operatorname{\sf Teich}}
\newcommand{\Hyp}{\operatorname{\sf Hyp}}

\newcommand{\Comp}{\operatorname{\sf Comp}}
\newcommand{\St}{\operatorname{St}}
\newcommand{\Per}{\operatorname{\sf Per}}
\newcommand{\Perspace}{\operatorname{{\Bbb P}\sf er}}

\newcommand{\End}{\operatorname{End}}

\newcommand{\Gr}{\operatorname{Gr}}
\newcommand{\Stab}{\operatorname{Stab}}

\newcommand{\Id}{\operatorname{Id}}

\newcommand{\Pic}{\operatorname{Pic}}
\newcommand{\Pos}{\operatorname{Pos}}

\newcommand{\Aut}{\operatorname{Aut}}

\newcommand{\Hdg}{{\operatorname{\sf Hdg}}}

\newcommand{\Diff}{\operatorname{\sf Diff}}

\newcommand{\NS}{\operatorname{NS}}

\renewcommand{\Re}{\operatorname{Re}}
\renewcommand{\Im}{\operatorname{Im}}

\newcounter{Mycounter}[section]
\newcounter{lemma}[section]
\newcounter{claim}[section]
\newcounter{sublemma}[section]
\newcounter{corollary}[section]
\renewcommand{\thecorollary}{{Corollary \thesection.\arabic{corollary}}}
\newcommand{\corollary}{%
    \setcounter{corollary}{\value{Mycounter}}
    \refstepcounter{corollary}
    \stepcounter{Mycounter}
    {\noindent \bf \thecorollary:\ }}

\newcounter{theorem}[section]
\renewcommand{\thetheorem}{{Theorem \thesection.\arabic{theorem}}}
\newcommand{\theorem}{%
    \setcounter{theorem}{\value{Mycounter}}
    \refstepcounter{theorem}
    \stepcounter{Mycounter}
    {\noindent \bf \thetheorem:\ }}

\newcounter{conjecture}[section]
\newcounter{proposition}[section]
\renewcommand{\theproposition}
      {{Proposition \thesection.\arabic{proposition}}}
\newcommand{\proposition}{%
    \setcounter{proposition}{\value{Mycounter}}
    \refstepcounter{proposition}
    \stepcounter{Mycounter}
    {\noindent \bf \theproposition:\ }}

\newcounter{definition}[section]
\renewcommand{\thedefinition}
      {{Definition~\thesection.\arabic{definition}}}
\newcommand{\definition}{%
    \setcounter{definition}{\value{Mycounter}}
    \refstepcounter{definition}
    \stepcounter{Mycounter}
    {\noindent \bf \thedefinition:\ }}

\newcounter{example}[section]
\newcounter{remark}[section]
\renewcommand{\theremark}{{Remark \thesection.\arabic{remark}}}
\newcommand{\remark}{%
    \setcounter{remark}{\value{Mycounter}}
    \refstepcounter{remark}
    \stepcounter{Mycounter}
    {\noindent \bf \theremark:\ }}

\newcounter{problem}[section]
\newcounter{question}[section]
\makeatletter

\@addtoreset{equation}{section} \@addtoreset{footnote}{section}
\makeatother

\def\blacksquare{\hbox{\vrule width 5pt height 5pt depth 0pt}}
\def\endproof{\blacksquare}

\begin{document}

\begin{center}
{\Large\bf
Collections of parabolic orbits\\[2mm] in homogeneous spaces, 
homogeneous dynamics \\[2mm] and hyperk\"ahler geometry\\[2mm]
}

Ekaterina Amerik\footnote{Partially supported by 
RSCF, grant number 14-21-00053.},  
Misha
Verbitsky\footnote{Partially supported by RSCF, grant number 14-21-00053.


}

\end{center}

{\small \hspace{0.15\linewidth}
\begin{minipage}[t]{0.7\linewidth}
{\bf Abstract} \\
Consider the space $M=O(p,q)/O(p)\times O(q)$
of positive $p$-dimensional subspaces in a pseudo-Euclidean space $V$ of signature $(p,q)$, where $p>0$, $q>1$ and $(p,q)\neq (1,2)$, with integral structure: $V=V_{\Z}\otimes \R$.
Let $\Gamma$ be an arithmetic subgroup in $G=O(V_{\Z})$, and 
$R\subset V_{\Z}$ a $\Gamma$-invariant
set of vectors with negative square. Denote by $R^\bot$ the 
set of all positive $p$-planes $W\subset V$ such that
the orthogonal complement $W^\bot$ contains $r\in R$.
We prove that either $R^\bot$ is dense in $M$ or $\Gamma$
acts on $R$ with finitely many orbits.
This is used to prove that the squares of primitive classes giving the rational 
boundary of the
K\"ahler cone (i.e. the classes of ``negative'' minimal rational curves) on a hyperkahler manifold $X$ are bounded by a number which
depends only on the deformation class of $X$.
We also state and prove the density of orbits in a more general 
situation when $M$ is the space of maximal compact subgroups in 
a simple real Lie group.
\end{minipage}
}

\tableofcontents


\section{Introduction}



In \cite{_AV:Mor_Kaw_},
the following result of hyperbolic geometry was established.

\hfill

\theorem \label{_density_hyperbo_Theorem_}
Let $\H=SO(1,n)/SO(n)$ be the 
$n$-dimensional real hyperbolic space, $n >2$, and 
$\Gamma\subset SO(1,n)$ an arithmetic subgroup naturally
acting on $\H$. 
Consider a $\Gamma$-invariant
collection $\{S_i\}$ of rational hyperplanes 
$S_i\subset \H$, and let $Z:= \bigcup_i S_i$
be their union in $\H$. Assume that $\Gamma$
acts on $\{S_i\}$ with infintely many orbits.
Then $Z= \bigcup_i S_i$ is dense in $\H$.

\hfill

Here we view $\H$ as the projectivization of the positive cone $V^+$ in the vector space $V$ with integral structure,
$V=\Z^{n+1}\otimes \R$, equipped with an integral quadratic
form of signature $(+, -, \dots, -)$).

This result can be understood as a statement about homogeneous
geometry of the space $M=G/K$,
where $G=SO(1,n)$ and $K$ a maximal compact subgroup; the
hyperplanes $S_i$ are orbits, under maximal parabolic
subgroups $P_i=\Stab(v_i)$, $v_i^2<0$,
of certain points $x_i\in v_i^{\bot}$.
In the present paper we generalize this result for
$SO^+(p,q)$ (connected component of $SO(p,q)$) for all $p, q \geq 2$.

\hfill

\theorem\label{_main_intro_p,q_Theorem_}
Let $G=SO^+(V)$ be a connected component
of the orthogonal group, where $V=V_\Z\otimes_\Z \R$ is a pseudo-Euclidean
space of  signature $(p,q)$, where $p>0$, $q>1$ and $(p,q)\neq (1,2)$, obtained from an integral
lattice $V_\Z$, $\Gamma\subset SO^+(V_\Z)$ an arithmetic
subgroup of $G$, and $R\subset V_\Z$ a $\Gamma$-invariant
set of vectors with negative squares. Denote by
$\Gr_+=O(p,q)/O(p)\times O(q)$
the space of all positive $p$-dimensional planes in $V$,
and let $R^\bot$ be the set of all planes $W\in \Gr_+$
such that the orthogonal complement $W^\bot$ contains
some $r\in R$. Then either $\Gamma$ acts on $R$ with 
finitely many orbits, or $R^\bot$ is dense on $\Gr_+(V)$.

{\bf Proof:} It is a special case of
\ref{_density_P_i_general_Theorem_}. \endproof

\hfill

Below we state and prove a generalization of this result
when $G$ is any simple non-compact algebraic Lie group. For our
applications, only the case $G=SO^+(3,q), q>1$ is
interesting, but it seems that this generatization
simplifies the problem and makes it conceptually easier.

\hfill
  
Our motivation comes from
the algebraic geometry of hyperk\"ahler manifolds; see Section \ref{HK}
for details on those.
The connection is as follows: on the second cohomology of a hyperk\"ahler manifold $X$, there is an integral non-degenerate
quadratic form $q$, the {\bf Beauville-Bogomolov-Fujiki (BBF) form}. The signature of $q$ is $(3, b_2-3)$. If $M$ is projective,
$q$ is of signature $(+,-,\dots,-)$ on the real Neron-Severi group $\NS(X)\otimes \R$; we view the projectivization of the positive cone
in $\NS(X)\otimes \R$ as the hyperbolic space $\H$. An important question of algebraic geometry is to describe the ample cone inside
the positive cone. It is well-known that it is cut out by a (possibly infinite) number of rational hyperplanes. One may ask whether there are
only finitely many of them up to the action of automorphism group (this is a version of the ``cone conjecture'' by Kawamata and 
Morrison). Our theorem from \cite{_AV:Mor_Kaw_} implies this as soon as the Picard
rank is greater than three. The group $\Gamma$ of the theorem is the Hodge 
monodromy group (see subsection \ref{monodr})  rather than
the automorphism group, but this is handled using
the global Torelli theorem (\cite{_Markman:survey_}, 
\cite{_V:Torelli_}).

Our motivation for the present paper is a refinement of the cone conjecture in the topological context. In fact we have shown 
in \cite{_AV:MBM_} that the K\"ahler cone inside the positive cone is a connected component of the complement to the union of hyperplanes 
orthogonal to the so-called {\bf MBM classes} of type $(1,1)$; those MBM classes are simply the classes whose orthogonal hyperplane
supports a face of the K\"ahler cone of a birational model of $X$, as well as their monodromy transforms (so that the set $R$ of MBM classes
is invariant by $\Gamma$). In this case, \ref{_density_hyperbo_Theorem_} says that the set of primitive 
MBM classes of type $(1,1)$ is finite up to the action of the monodromy group; since the latter acts by isometries, a consequence of \ref{_density_hyperbo_Theorem_} is
that the BBF squares of primitive MBM classes of type $(1,1)$  are bounded in absolute value by a constant $N$, which apriori depends on the complex 
manifold $M$. 

But we also have shown that the MBM property is deformation invariant, that is,
an MBM class remains MBM on all deformations where it stays of Hodge type $(1,1)$. It therefore makes sense to introduce the notion of an MBM class
in $H^2(X,\Z)$ without specifying its Hodge type, by requiring it to be MBM on those deformations where it is of type $(1,1)$. One then can ask
whether there exists an upper bound $N$ for the absolute value of
the BBF square of a primitive MBM class
which depends not on the complex structure, but only on its deformation type. By a result of Huybrechts, this is the same as to
ask that the bound only depends on topology (indeed the result affirms that there are only finitely many deformation types in the 
topological class). We have conjectured the affirmative answer in \cite{_AV:MBM_} (Conjecture 6.4).

The purpose of this paper is to prove this conjecture. The main ingredient of the proof
is the generalization of \ref{_density_hyperbo_Theorem_} to the space
$M=G/K$, where 
$G=SO(3,n)$ and $K= SO(3)\times SO(n)$. This is exactly \ref{_main_intro_p,q_Theorem_}
with $p=3$ which we state below separately in order to introduce some notations and 
terminology.

Recall that {\bf a lattice} in a Lie group
is a discrete subgroup of finite covolume
(that is, the quotient has finite volume). Arithmetic subgroups of reductive groups without non-trivial rational 
characters are lattices by Borel and Harish-Chandra theorem. 


\hfill

\theorem\label{_density_Gr_+++_Theorem_}
Let $\Gamma$ be an arithmetic lattice in $G$, where $G=SO^{+}(3, n)$, $n\geq 2$,
is the connected component of the unity of the group of linear isometries of a vector space
$(V,q)$ of signature $(3,n)$. Consider
a $\Gamma$-invariant set of rational vectors $R\subset V$
with negative square. Let 
$\Gr_{+++} =\frac{SO(3,n)}{SO(3)\times SO(n)}$
be the Grassmannian of 3-dimensional positive
oriented planes in $V$. For each $r\in R$,
denote by $S_r$ the set of all 
3-planes $W\in \Gr_{+++}$ orthogonal to $r$.
Assume that $\Gamma$ acts on $R$ with 
infinitely many orbits. Then the union
$\bigcup_{r\in R} S_r$ is dense in 
$\Gr_{+++}$.

\hfill

A more general theorem is proved in Section \ref{groups}, and the application 
to hyperk\"ahler geometry is explained in detail in Section \ref{HK}.
Here we just state the main corollary and give the idea of its proof.

\hfill

\corollary\label{MBM-bound} Let $X$ be a hyperk\"ahler manifold with $b_2(X)\geq 5$. The monodromy group acts with finitely many orbits on the set of primitive MBM classes in $H^2(X,\Z)$. The BBF square
of a primitive MBM class on $M$ is therefore bounded in absolute value by a constant $N$ which depends only on the topology of $X$.

\hfill

{\bf Idea of proof:} Let $R$ be the set of primitive MBM classes; then the union $\bigcup_{r\in R} S_r$ cannot be dense in $\Gr_{+++}$ since its
complement is identified to a connected component of the {\bf Teichm\"uller space of hyperk\"ahler structures}, open in $\Gr_{+++}$ (\cite{_AV:Teich_}).
\endproof

\hfill

\remark This corollary also removes a technical assumption $b_2\neq 5$ needed in \cite{_AV:Mor_Kaw_} to prove that the set of faces of the K\"ahler cone is finite up to automorphism group action; see section \ref{HK}.

\hfill

Notice that each $S_r$ is an orbit of some point $x_r$ in
$\Gr_{+++}$  under the maximal parabolic subgroup $P_r=\Stab(r)$; such an orbit is
special in the sense that the subspace corresponding to $x_r$ is orthogonal to $r$.
By analogy with the hyperbolic space setting, we call them
{\bf parabolic orbits of hyperplane type}.

\medskip

\ref{_density_Gr_+++_Theorem_} and 
\ref{_density_hyperbo_Theorem_}
are special cases of a more general statement.
Denote by $G$ a connected simple real algebraic Lie group, $\g$ its
Lie algebra, and $\g_\C$ its complexification.

Recall that a connected subgroup $P\subset G$
is called {\bf parabolic}
if the complexification ${\goth p}_\C$ of its
Lie algebra contains a Borel subalgebra of
$\g_\C$, and {\bf maximal parabolic} when it is maximal with these properties. In particular there are no other connected 
closed subgroups between $P$ and $G$.

Let $K$ be a maximal compact subgroup of $G$.
By Cartan's theorem, $G$ retracts on $K$, so
$K$ is itself connected (\cite[Theorem 3.1, Chapter XV]{_Hochschild:Lie_}).
Since the normalizer of a maximal compact
subgroup of a simple algebraic Lie group
is again compact, $K$ is equal to its
normalizer, as follows from the Iwasawa decomposition. 
Therefore, we have a natural bijection between $G/K$ and
the set of conjugates of $K$, that is, one may view $G/K$ as a set of 
all maximal compact subgroups of $G$. The action of $G$ on $G/K$ by 
left translation becomes, under this identification, the action of $G$ on the 
set of maximal compact subgroups by conjugation.

Let $P$ be a maximal parabolic subgroup and $K'$ a maximal compact. Call $P$ and $K'$ {\bf compatible} if $P\cap K'$ is maximal compact in $P$. If $P$ and $K'$ are 
compatible, then clearly $P$ is compatible with all maximal compacts in the 
$P$-orbit of $K'$. 

\hfill

\definition\label{hpl-type}
An {\bf orbit of hyperplane type} under the action of $P$ on $G/K$ is a $P$-orbit consisting of $P$-compatible maximal
subgroups.  

\hfill

Now we can formulate the general statement.

\hfill

\theorem\label{_density_P_i_general_Theorem_}
Let $\Gamma$ be an arithmetic lattice in $G$, where $G$
is a connected simple real algebraic Lie group, $K$ its maximal compact subgroup,
and $P$ a maximal parabolic subgroup which is assumed to
be generated by unipotents. We assume that the groups $G,
K, P$ are defined over $\Q$.
Consider a $\Gamma$-invariant set $\{S_i\}$ of orbits of hyperplane type of subgroups 
$x_i P x_{i}^{-1}$ acting on $G/K$.
Assume that $\Gamma$ acts on $\{S_i\}$ with 
infinitely many orbits. Then the union
$\bigcup S_r$ is dense in 
$G/K$.

\hfill

The next section is devoted to the proof of \ref{_density_P_i_general_Theorem_}.



\section{Homogeneous dynamics, parabolic subgroups and Mozes-Shah theorem}\label{groups}



We deduce \ref{_density_P_i_general_Theorem_}
from the general formalism of Ratner, Mozes-Shah and Eskin-Mozes-Shah
used in \cite{_AV:Mor_Kaw_}  to prove \ref{_density_hyperbo_Theorem_}. 

In order to be able to apply this machinery we first prove  
a simple statement on Lie groups which replaces a set of orbits of hyperplane type of conjugated parabolic 
subgroups on $G/K$ by a set of orbits of a single one on a suitable fibration over $G/K$, which is a $G$-homogeneous
space ``in between'' $G/K$ and $G$ itself.

As above, we identify $G/K$ with the set of all maximal
compact subgroups of $G$.

\hfill

\proposition\label{_parabo_orbi_Proposition_}
Let $G$ and $K$ be as above. Consider a maximal
parabolic subgroup $P\subset G$, and let $S_i=y_iPy_i^{-1} x_i$,
$i\in I$  by a set of orbits of hyperplane type under 
conjugates of $P$. Then there exists a single $P_0$
conjugate to $P$ such that the $S_i$ are projections of
$P_0$-orbits $R_i$ in $G$. 

\hfill

\remark Those projections themselves are of course not $P_0$-orbits, 
as the projection map we 
consider is not equivariant. In the case when $G$ is $SO^+(1,2)$, this construction
is known as the {\bf geodesic flow}: the hyperplanes in the hyperbolic plane lift 
tautologically to the unit tangent bundle as orbits of a single $SO^+(1,1)$, 
corresponding to the subgroup of diagonal matrices under the identification with
$PSL(2,\R)$. 

\hfill

{\bf Proof of \ref{_parabo_orbi_Proposition_}:} Under the identifications we have made, the action 
of $G$ on $G/K$ corresponds to the adjoint action on $M$, that is, $x\in G$
sends a maximal compact subgroup $K$ to $xKx^{-1}$. Consider the space
$M_1$ of all pairs $(K_1, P_1)$, where $K_1\subset G$ is a maximal 
compact subgroup, $P_1\subset G$ a maximal parabolic subgroup
which is conjugate to $P$, and $K_1\cap P_1$ a maximal compact subgroup
of $P_1$; that is, $K_1$ and $P_1$ are compatible.  

(For instance, in the situation of \ref{_density_Gr_+++_Theorem_},
$G=SO(3,n)$, $K=SO(3)\times SO(n)$, the maximal parabolic group
$P_1=\St_G(v)$ is a stabilizer of a vector $v$ with negative square,
and $K_1\cap P_1$ is isomorphic to $SO(3)\times SO(n-1)$.) 

Now let $M_0$ denote the space of maximal parabolic
subgroups conjugate to $P$, that is, subgroups of the
form $P_1 = gPg^{-1}\subset G$. There is a natural diagram
\begin{diagram}[size=2em]
& & M_1 &&\\
&\ldTo^{\pi} && \rdTo^{\sigma}&\\
M &&&& M_0
\end{diagram}
with forgetful maps $\pi$ and $\sigma$.
Since the normalizer $N(P)$ is an algebraic group which 
has the same connected component of the unity as $P$, and $M_0$ is naturally
identified with $G/N(P)$, the standard map $G/P\arrow M_0=G/N(P)$
is a finite covering. Therefore the orbits of hyperplane type
under $P_1 = gPg^{-1}$ on $M$ are connected
components of $\pi(\sigma^{-1}(P_1))$. As $P_1$ varies in $M_0$, the 
connected components of $\sigma^{-1}(P_1)$ give 
a $G$-invariant foliation on $M_1$ with leaves which
are mapped to $gPg^{-1}$-orbits of hyperplane type in $M$. Taking the preimages of those leaves in $G$, 
we obtain a translation-invariant foliation on $G$ with leaves which
are mapped to $gPg^{-1}$-orbits in $M$. But such a foliation is necessarily
by orbits of the action of a subgroup; this subgroup is the $P_0$ that we are looking for.
\endproof

\hfill

{\bf Proof of \ref{_density_P_i_general_Theorem_}:}\\
We have seen that the collection of $S_i$
lifts to a $\Gamma$-invariant collection of $P_0$-orbits $R_i$ in $G$. It suffices to prove that this collection is dense or finite up to
the action of $\Gamma$, or, in other words, that the corresponding orbits in $\Gamma\backslash G$ are finitely many or dense.
To this end, we apply the same argument as in
\cite{_AV:Mor_Kaw_}.
The $P_0$-orbits $R_i$ give rise to probability 
measures $\mu_i$ on 
$\Gamma\backslash G$, supported on those orbits; 
those are simply the translates of $\mu_{P_0}$, the ``pushforward'' of the
Haar measure on $P_0$. If there are infinitely many 
of them, then by Mozes-Shah and Dani-Margulis theorems
(\cite{_Mozes_Shah_}, Corollaries 1.1, 1.3, 1.4)  one can
extract from $\{\mu_i\}$ a weakly converging subsequence of
measures on the one-point compactification of
$\Gamma\backslash G$,  possibly converging
to the measure $\mu_{\infty}$, concentrated at the infinite point,
but otherwise to a probability measure on an orbit of another closed subgroup
$P_0'$ containing $P$ or its conjugate (\cite{_Mozes_Shah_}, Theorem 1.1).

We remark that the results of Mozes-Shah from
\cite{_Mozes_Shah_} are valid for any collection of 
measures ergodic with respect to subgroups generated by unipotents (as can be found in the same paper, see e.g. Lemma 2.3). 
The crucial point for our theorem is that our converging subsequence
consists of measures supported on orbits of the same
maximal parabolic 
subgroup $P_0$, which is of finite index in its normalizer. 
Then it cannot converge to infinity 
by a theorem of Eskin-Mozes-Shah (\cite{EMS}, Theorem 1.1). Convergence to another translate of $\mu_{P_0}$
is impossible because of the finiteness of the index of  $P_0$ in its normalizer, as in Lemma 4.12 of \cite{_AV:Mor_Kaw_}.
 Therefore it converges to a translate of the ergodic measure 
$\mu_Q$ 
associated to a closed subgroup $Q$ strictly
containing (a conjugate of) $P$; but such $Q$ can only be $G$ itself, from where the density.

This finishes the proof of  \ref{_density_P_i_general_Theorem_}.



\section{Teichm\"uller space for hyperk\"ahler structures 
and MBM classes}\label{HK}


\subsection{Hyperk\"ahler manifolds and monodromy}\label{Monodr}


In this subsection, we recall some basic results on hyperk\"ahler manifolds.

\hfill

\definition
A {\bf hyperk\"ahler manifold}
is a compact K\"ahler holomorphically symplectic manifold.

\hfill

\definition
A hyperk\"ahler manifold $M$ is called
{\bf simple}, or IHS, if $\pi_1(M)=0$, $H^{2,0}(M)=\C$.

\hfill

This definition is motivated by Bogomolov's decomposition theorem:

\hfill

\theorem \label{_Bogo_deco_Theorem_}
(\cite{_Bogomolov:decompo_})
Any hyperk\"ahler manifold admits a finite covering
which is a product of a torus and several 
simple hyperk\"ahler manifolds.
\endproof

\hfill



\remark
Further on, we shall assume that
all hyperk\"ahler manifolds we consider are 
of maximal holonomy, that is, simple.

\hfill

An important property of hyperk\"ahler manifolds is the presence of an integral quadratic form of signature $(3,b_2-3)$ 
on their second
cohomology, the {\bf Bogomolov-Beauville-Fujiki (BBF) form}. It was
defined in \cite{_Bogomolov:defo_} and 
\cite{_Beauville_} using integration of differential forms,
but it is easiest to describe it by the
Fujiki theorem, proved in \cite{_Fujiki:HK_}; it stresses the topological origin of the BBF form.

\hfill

\theorem\label{_Fujiki_Theorem_}
(Fujiki)
Let $M$ be a simple hyperk\"ahler manifold,
$\eta\in H^2(M)$, and $n=\frac 1 2 \dim M$. 
Then $\int_M \eta^{2n}=c q(\eta,\eta)^n$,
where $q$ is a primitive integral quadratic form on $H^2(M,\Z)$,
and $c>0$ a constant (depending on $M$). \endproof

\hfill

The signature of the BBF form on the space of $(1,1)$-classes is $(+,-,\dots, -)$.

\hfill



\definition\label{_posi_cone_Definition_}
A cohomology class $\eta \in H^2_\R(M)$ is called
{\bf negative} if $q(\eta,\eta)<0$, and {\bf positive}
if $q(\eta,\eta)>0$. The {\bf positive cone} $\Pos(M)\in H^{1,1}(M)$ is 
one of the two connected components of the set of positive $(1,1)$-classes which contains the K\"ahler classes.




\hfill

\definition
Let $M$ be a compact complex manifold, and 
$\Diff_0(M)$ a connected component of its diffeomorphism group
({\bf the group of isotopies}). Denote by $\Comp$
the space of complex structures of K\"ahler type on $M$ (remark here that the set of complex structures of K\"ahler type 
is open in the space of all complex structures by Kodaira-Spencer stability theorem), and let
$\Teich:=\Comp/\Diff_0(M)$. We call 
it {\bf the Teichm\"uller space.} 

\hfill

For hyperk\"ahler manifolds, this is a finite-dimensional complex
non-Hausdorff manifold (\cite{_Catanese:moduli_}, \cite{_V:Torelli_}).

\hfill

\definition The {\bf mapping class group} is 
$\Diff(M)/\Diff_0(M)$, naturally acting  on $\Teich$. 


\hfill

It follows from a result of Huybrechts (see \cite{_Huybrechts:finiteness_}) 
that in the hyperk\"ahler case $\Teich$ has only finitely many
connected components. Therefore, the subgroup of the mapping class
group which fixes the connected component of our chosen complex structure
is of finite index in the mapping class group. 

\hfill

\definition\label{monodr} The {\bf monodromy group}
$\Gamma$ is the image of this subgroup in $\Aut H^2(M,\Z)$. 
The {\bf Hodge monodromy group}
is the subgroup $\Gamma_\Hdg\subset \Gamma$ preserving the Hodge decomposition.

\hfill

\theorem\label{arithmetic} (\cite{_V:Torelli_}, Theorem 3.5) 
The monodromy group is a finite index subgroup in $O(H^2(M, \Z), q)$
(and the Hodge monodromy therefore acts as an arithmetic subgroup of the orthogonal group on the 
Picard lattice).

\hfill

In \cite{_Markman:survey_},,
E. Markman obtained the following
implication of the global Torelli theorem relating 
the Hodge monodromy to automorphisms.

\hfill

\theorem\label{markman-torelli} If $\gamma\in \Gamma_\Hdg$ takes a K\"ahler class to a K\"ahler class, then
$\gamma=f^*$ for some $f\in \Aut(M)$.

\subsection{MBM classes}

The notion of an MBM class was introduced in
\cite{_AV:MBM_} 
in order to understand better how the K\"ahler cone 
sits in the positive cone.

\hfill

\definition
An integral $(1,1)$-class $z$ on $M$ is MBM if for some $\gamma\in \Gamma_\Hdg$, the hyperplane 
$\gamma(z)^{\bot}$ supports a (maximal-dimensional) face of the K\"ahler cone of a birational model of $M$.

\hfill

This somewhat mysterious definition has a simple geometric
interpretation thanks to the deformation invariance
of the MBM property. Note that the BBF form identifies
homology with rational coefficients and cohomology.
The following two theorems are proven in \cite{_AV:MBM_}.

\hfill

\theorem Let $M$ be a non-algebraic hyperk\"ahler manifold
with $\Pic(M)=\langle z \rangle$, where $z\in H_2(M,\Z)$ a
negative homology class. Then $z$ is an MBM class
if and only if for some $\lambda\neq 0$, $\lambda z$ can
be represented by a curve. 

\hfill

\theorem
Let $M, M'$ be hyperk\"ahler manifolds in the same
deformation class, such that a negative cohomology class $z\in
H^2(M, \Z)$ is of type (1,1) on $M$ and $M'$. Then
$z$ is an MBM class on $M$ if and only if it is MBM on
$M'$.

\hfill

This result allows one to extend the notion of MBM classes to
the whole of $H^2(M,\Z)$.








\hfill

\definition
A negative class $\eta\in H^2(M,\Z)$ is called {\bf MBM}
if it is MBM for some deformation of $M$ where it is of type $(1,1)$.

\hfill

By the very definition, the MBM property is deformation-invariant.

\hfill

\remark It follows immediately from the definition that the set of MBM classes in this generalized sense 
in $H^2(M, \Z)$ is also $\Gamma$-invariant, where $\Gamma$ is the 
monodromy group of $M$. Indeed given an MBM class $z$ and a monodromy transform $\gamma(z)$ one can find a deformation of $M$ such that both
of them are of type $(1,1)$; then $\gamma(z)$ is MBM in the sense of our first definition.

\hfill

The main result of this paper is the following theorem.

\hfill

\theorem\label{_finite_MBM_Theorem_}
Let $R\subset H^2(M,\Z)$ be the set of primitive MBM classes
in the cohomology of a hyperk\"ahler manifold whose second
Betti number $b_2(M)$ is at least 5. Then the
monodromy group $\Gamma$ acts on $R$ with finitely many
orbits. In particular, there is a number $N$ depending
only on the deformation type of $M$, such that for any
$z\in R$, $|q(z)|\leq N$.

\hfill
 
{\bf Proof:} Consider $G=SO^+(3,n)$ where $n=b_2(M)-3>1$,
and the homogeneous space $G/K$ as in 
\ref{_density_Gr_+++_Theorem_}, identified with the
grassmannian of positive 3-planes in the vector 
space  $H^2(M,\R)$ of signature $(3,n)$.
The subsets formed by 3-planes orthogonal to a given MBM
class in $H^2(M,\Z)$ form an orbit of hyperplane type. As $n>1$, by
\ref{_density_Gr_+++_Theorem_} either the set of MBM
classes is finite up to $\Gamma$-action, or the corresponding
orbits of hyperplane type are dense in $G/K$. The latter
is impossible. Indeed, in the
same way as the orthogonals to the MBM classes of type
$(1,1)$ serve as boundaries of the K\"ahler cone
inside the positive cone, the orthogonals to the all
MBM classes are a complement of a meaningful object 
of modular nature, the Teichm\"uller space of 
hyperk\"ahler structures (\cite{_AV:Teich_}), open in the
Grassmannian. In the next and final subsection, we briefly
recall this for reader's convenience.



\subsection{Hyperk\"ahler Teichm\"uller space}
\label{_hype_Teich_Subsection_}

\definition
Let $(M,g)$ be a Riemannian manifold, and $I,J,K$
endomorphisms of the tangent bundle $TM$ satisfying the
quaternionic relations
\[
I^2=J^2=K^2=IJK=-\Id_{TM}.
\]
The triple $(I,J,K)$ together with
the metric $g$ is called {\bf a hyperk\"ahler structure}
if $I, J$ and $K$ are integrable and K\"ahler with respect to $g$.

Consider the K\"ahler forms $\omega_I, \omega_J, \omega_K$
on $M$:
\[
\omega_I(\cdot, \cdot):= g(\cdot, I\cdot), \ \
\omega_J(\cdot, \cdot):= g(\cdot, J\cdot), \ \
\omega_K(\cdot, \cdot):= g(\cdot, K\cdot).
\]
An elementary linear-algebraic calculation implies
that the 2-form $\Omega:=\omega_J+\1\omega_K$ is of Hodge type $(2,0)$
on $(M,I)$. This form is clearly closed and
non-degenerate, hence it is a holomorphic symplectic form.

\hfill

In algebraic geometry, the word ``hyperk\"ahler''
is essentially synonymous with ``holomorphically
symplectic'', due to the following theorem, which is
implied by Yau's solution of Calabi conjecture
(\cite{_Besse:Einst_Manifo_}, \cite{_Beauville_}).

\hfill

\theorem\label{_Calabi-Yau_Theorem_}
Let $M$ be a compact, K\"ahler, holomorphically
symplectic manifold, $\omega$ its K\"ahler form, $\dim_\C M =2n$.
Denote by $\Omega$ the holomorphic symplectic form on $M$.
Suppose that $\int_M \omega^{2n}=\int_M (\Re\Omega)^{2n}$.
Then there exists a unique hyperk\"ahler metric $g$ with the same
K\"ahler class as $\omega$, and a unique hyperk\"ahler structure
$(I,J,K,g)$, with $\omega_J = \Re\Omega$, $\omega_K = \Im\Omega$.
\endproof

\hfill

Every hyperk\"ahler structure induces a whole 2-dimensional
sphere of complex structures on $M$, as follows. 
Consider a triple $a, b, c\in R$, $a^2 + b^2+ c^2=1$,
and let $L:= aI + bJ +cK$ be the corresponging quaternion. 
Quaternionic relations imply immediately that $L^2=-1$,
hence $L$ is an almost complex structure. 
Since $I, J, K$ are K\"ahler, they are parallel with respect
to the Levi-Civita connection. Therefore, $L$ is also parallel.
Any parallel complex structure is integrable, and K\"ahler.
We call such a complex structure $L= aI + bJ +cK$
{\bf a complex structure induced by the hyperk\"ahler structure}.
There is a 2-dimensional holomorphic family of 
induced complex structures, and the total space
of this family is called {\bf the twistor space}
of a hyperk\"ahler manifold, its base being {\bf the 
twistor line} in the Teichm\"uller space
$\Teich$ which we are going to define next.

\hfill

\definition
Let $(M,I,J,K,g)$ and $(M,I',J',K',g')$ be two
hyperk\"ahler structures. We say that these structures
are {\bf equivalent} if the corresponding quaternionic algebras
in $\End(TM)$ coincide.

\hfill

Consider the infinite-dimensional space $\Hyp$ of all quaternionic
triples $I,J,K$ on $M$ which are induced by some
hyperk\"ahler structure, with the same $C^\infty$-topology
of convergence with all derivatives. The quotient 
$\Hyp/SU(2)$ (which is probably better to write as
$\Hyp/SO(3)$, since $-1$ acts trivially
on the triples)  is naturally identified with the set
of equivalence classes of hyperk\"ahler structures, 
up to changing the metric $g$ by a constant multiplier.

\hfill

\remark\label{_metric_fixed_vol_hk_Remark_}
As shown in \cite{_AV:Teich_}, for hyperk\"ahler manifolds
with maximal holonomy the quotient $\Hyp_m:=\Hyp/SU(2)$  is
also identified with the space of  all hyperk\"ahler
metrics of fixed volume, say, volume 1.

\hfill

\definition\label{def-teich-hk}
Define {\bf the Teichm\"uller space of 
hyperk\"ahler structures} as the quotient $\Hyp_m/\Diff_0$,
where $\Diff_0$ is the connected component of the group
of diffeomorphisms $\Diff$, and {\bf the moduli
of hyperk\"ahler structures} as $\Hyp_m/\Diff$.

\hfill

\definition\label{_Teich_h_Definition_}
Let $M$ be a hyperk\"ahler manifold of maximal holonomy,
and  $\Teich_h:= \Hyp_m/\Diff_0$ the Teichm\"uller
space of hyperk\"ahler structures. Consider the
space $\Perspace_h=Gr_{+++}(H^2(M,\R))$ of all positive 
oriented 3-dimensional subspaces
in $H^2(M,\R)$, naturally diffeomorphic to
$\Perspace_h\cong SO(b_2-3, 3)/SO(3)\times
SO(b_2-3)$. Let $\Per_h:\; \Teich_h \arrow\Perspace_h$ 
be the map associating  the 3-dimensional space generated
by the three K\"ahler forms $\omega_I, \omega_J,\omega_K$
to a hyperk\"ahler structure $(M,I,J,K,g)$.
This map called {\bf the period map for the 
Teichm\"uller space of hyperk\"ahler structures},
and $\Perspace_h$ {\bf the period space of  hyperk\"ahler
  structures}.

\hfill

\theorem\label{th-teich-hk}
Let $M$ be a hyperk\"ahler manifold of maximal holonomy,
and $\Per_h:\; \Teich_h \arrow\Perspace_h$ the period map
for the 
Teichm\"uller space of hyperk\"ahler structures.
Then the period map 
$\Per_h:\; \Teich_h \arrow\Perspace_h$ is an open embedding
for each connected component. Moreover, its
image is the set of all spaces $W\in \Perspace_h$
such that the orthogonal complement $W^\bot$
contains no MBM classes.

{\bf Proof:} See \cite{_AV:Teich_}.
\endproof

\hfill

Let $V=H^2(M,\R)$ be the second cohomology of
a hyperk\"ahler manifold, equipped with the BBF form.
Denote by $\Gr_{+++}$ the space of all positive oriented 
3-planes in $V$. For each primitive MBM class $x$, the set of $W\in \Gr_{+++}$
orthogonal to $x$ is an orbit of the maximal parabolic 
group $P_x=\St_{SO(V)}(x)$. By \ref{th-teich-hk}, 
the space $\Teich_h$ is identified with $\Gr_{+++}\backslash \bigcup_{x\in R}P_x(W_x)$,
where $W_x\in \Gr_{+++}$ is any 3-space orthogonal to $x$, and $R$ the set
of primitive MBM classes. Denote by $\Gamma\subset SO(V)$
the monodromy group of $M$.
As shown in \cite{_V:Torelli_}, it is a lattice in $SO(V)$.
Since $P_x(W_x)$ is determined by $x$ and determines it 
uniquely, up to a sign, the number of $\Gamma$-orbits on $R$
is infinite if and only if the number of $\Gamma$-orbits on 
$P_x(W_x)$ is infinite. However, the union $\bigcup_{x\in R}P_x(W_x)$
is closed by \ref{th-teich-hk}, hence \ref{_density_Gr_+++_Theorem_} implies
that $R/\Gamma$ is finite. We have proved \ref{_finite_MBM_Theorem_}.

\hfill

{\bf Acknowledgements:}
We are grateful to Alex Eskin, Misha Kapovich and Maxim Kontsevich 
for illuminating discussions of hyperbolic geometry.

{

\noindent {\sc Ekaterina Amerik\\
{\sc Laboratory of Algebraic Geometry,\\
National Research University HSE,\\
Department of Mathematics, 7 Vavilova Str. Moscow, Russia,}\\
\tt  Ekaterina.Amerik@math.u-psud.fr}, also: \\
{\sc Universit\'e Paris-11,\\
Laboratoire de Math\'ematiques,\\
Campus d'Orsay, B\^atiment 425, 91405 Orsay, France}

\hfill

\noindent {\sc Misha Verbitsky\\
{\sc Laboratory of Algebraic Geometry,\\
National Research University HSE,\\
Department of Mathematics, 7 Vavilova Str. Moscow, Russia,}\\
\tt  verbit@mccme.ru}, also: \\
{\sc Universit\'e Libre de Bruxelles, CP 218,\\
Bd du Triomphe, 1050 Brussels, Belgium}
}

\end{document}